\documentclass[]{report}
\usepackage{amssymb}
\usepackage{titling}
\usepackage{amsmath}
\usepackage{geometry}
\usepackage{graphicx}
\usepackage{subfig}
\usepackage{fancyhdr}
\usepackage{lastpage}
\usepackage{lipsum}
\usepackage{ragged2e}
\usepackage{amsmath}
\usepackage[ruled,vlined]{algorithm2e}
\usepackage{algpseudocode}
\usepackage{comment}

\usepackage{hyperref}
\hypersetup{colorlinks,
	citecolor=blue
}
\usepackage{listings}
\usepackage{color}
\definecolor{DarkRed}{RGB}{139,0,0}
\definecolor{dkgreen}{rgb}{0,0.6,0}
\definecolor{gray}{rgb}{0.5,0.5,0.5}
\definecolor{mauve}{rgb}{0.58,0,0.82}

\usepackage{titling}

\title{\Large{Hybrid Dealiased Convolutions}}
\author{Robert Joseph\\ Department of Mathematical and Statistical Sciences}

\def\argmin{\mathop{\rm argmin}}
\begin{document}

\begin{titlepage}
    \begin{center}
        \vspace*{1cm}
            
        \Huge
        \textbf{Hybrid Dealiased Convolutions}
            
        \vspace{0.5cm}
        \LARGE
        An undergraduate thesis presented by            
        \vspace{0.5cm}  
        
        {Robert Joseph George}
        \vspace{0.5cm}
        
        Advised by: John Bowman and Noel Murasko\\
        Supervised by: Nicolas Guay
        \vspace{1cm}

        \normalsize{\hspace{12pt} \textbf{Abstract}}
        
        \justify 
        {This paper proposes a practical and efficient solution for computing convolutions using hybrid dealiasing. It offers an alternative to explicit or implicit dealiasing and includes an optimized hyperparameter tuning algorithm that uses experience to find the optimal parameters. Machine learning algorithms and efficient heuristics are also developed to estimate optimal parameters for larger convolution problems using only small squares/rectangles.}

        \vspace{1cm}
        \centering
        \normalsize{Submitted in partial fulfillment of the Honors requirements for the degree \\ of Bachelor of Honors in Applied Mathematics and Computer Science.}
            
        \vspace{0.8cm}
            
        \includegraphics[width=0.3\textwidth]{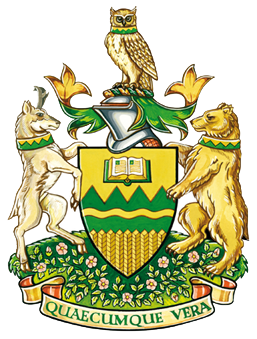}
        \vspace{0.5cm}

        \Large
        Department of Mathematics and Statistics \\
        University of Alberta\\
        \hspace{12pt} Edmonton, Canada
        \date{}
            
    \end{center}
\end{titlepage}
\hypersetup{linkcolor=DarkRed}
\tableofcontents
\pagebreak

\section{Introduction}
\vspace{-5pt}

Fast Fourier Transforms (FFT) and convolutions are important mathematical operations in many fields, including signal processing, image processing, computer vision and applied mathematics. FFTs are used to efficiently compute a signal's discrete Fourier transform (DFT), allowing for the analysis of the signal's frequency content. However, the convolution theorem allows for the efficient calculation of circular convolutions using the DFT. This theorem is useful in many areas of signal processing, such as image and audio processing, where convolution is used to apply filters and effects to signals. It also allows efficient convolution calculation in the frequency domain, which can be much faster than calculating a convolution in the time domain.

\justify
In most applications, we prefer computing a linear convolution instead of a circular convolution. Hence, using the convolution theorem, we pad the sequences with enough zeros to efficiently compute the linear convolution using the FFT algorithm. However, the caveat is that if the input sequences are not padded enough, then \textit{aliases} occur: errors that arise from the lack of periodicity in the input. To avoid aliases when using the convolution theorem in the context of linear convolution, we need to ensure that the input signals' frequency domain representation does not overlap. Aliasing occurs when different frequency components of signals overlap in the frequency domain, leading to distortions and errors in the resulting convolution.

\justify
One way to avoid aliasing is to ensure that the length of the output sequence after convolution is at least equal to the sum of the lengths of the input sequences minus one. This ensures that all the frequency components of the input signals are fully represented in the output sequence and there is no overlap between them. If necessary, we can zero-pad the input sequences before taking their Fourier transforms to ensure that the length of the output sequence meets this requirement. This approach is called \textit{Explicit dealiasing} and has become the standard practice.

\justify
Several approaches to avoid aliasing were discovered to bypass the explicit padding method, such as by using a different convolution method specifically designed to avoid overlap in the frequency domain, such as the overlap-add or overlap-save methods \cite{3}. These methods involve dividing the input signals into overlapping segments, performing convolution on each segment separately, and then combining the results using overlap-add or overlap-save techniques. This approach can be more computationally efficient than zero-padding and direct convolution, especially for long input signals.

\justify
An alternative to explicit dealiasing called \textit{Implicit dealiasing} \cite{2} \cite{3} \cite{4} can be used to account for known zero values in FFTs without the need for explicit zero padding. Previous work has focused on developing implicit dealiasing techniques for specific padding ratios, such as 1/2 and 2/3, for complex and Hermitian symmetric input data. These techniques are essential for applications such as signal processing and pseudo-spectral methods for solving partial differential equations. However, many other applications may not satisfy these specific requirements. Implicit dealiasing can be a more efficient approach for computing convolutions in the frequency domain, as it avoids the need for explicit zero-padding and the associated computational costs. By taking advantage of the known zero values in the input signals, implicit dealiasing techniques can help to improve the accuracy and efficiency of FFT-based convolution methods for a wide range of applications.

\justify
Although implicit dealiasing is an excellent approach, last summer \cite{1}, we developed and formulated a systematic framework that offered unique solutions by expanding \textit{Hybrid dealiasing} to issues with uneven input or minimum padded sizes, reducing memory and computation time. The work on equal cases (when both sequences/arrays have the same length) has been completed and is in submission \cite{2}. In this paper, we present the results of our work over the past year by expanding the hybrid dealiasing framework to incorporate the unequal case in $1/2/3$ dimensions, where we now have image convolution and construct the first hybrid dealiasing solution, which is much quicker than typical explicit padding approaches. We also expanded this to multi-convolution, applying the convolution to a sequence of data, $N$ arrays and demonstrating a conventional grid search hyperparameters (Algorithm $\textbf{A0}$) search strategy that selects the optimal hyperparameters that result in the shortest convolution time and is optimal. A general $N$ dimension approach was also developed, and we are developing them for existing state-of-the-art fast Fourier transform libraries such as FFTW++, which can result in low-memory FFT algorithms that are faster \footnote{Codebase will be released later in \href{https://github.com/dealias}{https://github.com/dealias}.}.

\justify
This semester, we worked on creating a more optimized and efficient hyperparameter tuning algorithm, which we call \textit{experience} to find the optimal parameters. This algorithm also searched a broader set of parameters to empirically determine the fastest algorithm for a given convolution problem, depending on the user. We also developed efficient heuristics, allowing us to estimate optimal parameters for more significant convolution problems using only small squares/rectangles.

\justify 
The remainder of this article is organized as follows: To understand the topic better, we begin with some preliminary background material and notations in Section 2. Section 3 introduces various dealiasing techniques, including our newly proposed Hybrid Dealiasing approach. Next, in Section 4, we discuss the optimization algorithms we have developed to improve the accuracy and efficiency of the dealiasing techniques. The methodology used in our research is presented in Section 5. Moving forward, Section 6 is dedicated to presenting the results of our research. Finally, in Sections 7 and 8, we discuss our work's limitations and potential applications, as well as future research directions.

\section{Fourier Transforms and  Convolutions }
\subsection{Discrete Fourier Transform}
\label{dft}
The Fourier transform (FT) of the function $f(x)$ is,
$
F(\omega)=\int_{-\infty}^{\infty} f(x) e^{-i \omega x} d x
$
and the inverse Fourier transform is
$
f(x)=\frac{1}{2 \pi} \int_{-\infty}^{\infty} F(\omega) e^{i \omega x} d \omega
$. When a signal is discrete and periodic, we consider the DFT. Let us consider a signal that has a period and whose length is $L$. Let us define the $M $th primitive root of unity as $\zeta_M=\exp \left(\frac{2 \pi i}{M}\right)$. We also now consider padding the input data to length $M$. Therefore, now we have a buffer $\boldsymbol{f} \doteq\left\{f_j\right\}_{j=0}^{M-1}$ where $f_j=a_j$ for $j<L$ and $f_j=0$ for all $j \geq L$. The discrete Fourier transform (DFT) of $f$ (forward transform) can be written as
\begin{equation}    
F_k = \frac{1}{M} \sum_{j=0}^{M-1} \zeta_M^{k j} f_j=\sum_{j=0}^{L-1} \zeta_M^{k j} f_j, \; k=0, \ldots, M-1. \
\end{equation} 

\justify
The corresponding backward DFT is given by \begin{equation}
    f_j = \frac{1}{M} \sum_{k=0}^{M-1} \zeta_M^{-k j} F_k, \; j=0, \ldots, M-1.
\end{equation}

\justify
It is easy to observe that this transform pair has an orthogonality condition, given by $$\sum_{j=0}^{M-1} \zeta_M^{\ell j}= \begin{cases}M & \text { if } \ell=s M \text { for } s \in \mathbb{Z} \\ \frac{1-\zeta_M^{\ell M}}{1-\zeta_M^{\ell}}=0 & \text { otherwise }\end{cases}$$

\subsection{Fast Fourier Transforms}
Hailed as one of the greatest algorithms of the 20th century, the FFT exploits the symmetry in the roots of unity by observing that $\zeta_M^r=\zeta_{M / r} \text { and } \zeta_M^M=1$. Intuitively, we break up the forward transform of length $M$ into two transforms of length $M/2$ using the identity $\forall k=0, \ldots, M-1$, which results in 
$$
F_k \doteq \frac{1}{M} \sum_{j=0}^{M-1} \zeta_M^{-k j} f_j = \frac{1}{M} \left(\sum_{j=0}^{\frac{M}{2}-1} \zeta_M^{-k 2j} f_{2j} + \sum_{j=0}^{\frac{M}{2}-1} \zeta_M^{-k (2j+1)} f_{2j+1}\right) = \frac{1}{M} \left(\sum_{j=0}^{\frac{M}{2}-1} \zeta_{M/2}^{-k j} f_{2j} + \zeta_M^{-k} \sum_{j=0}^{\frac{M}{2}-1} \zeta_{M/2}^{-kj} f_{2j+1}\right).
$$
This is extremely useful as it reduces the DFT's time complexity from $\mathcal{O}(M^2)$ to $\mathcal{O}(M\log M)$. For our case, since we care about convolutions, we need to pad the array with zeros as given by $\boldsymbol{f}$. We further assume that $L$ and $M$ share a common factor $m$; that is $L=p m$ and $M=q m$, where $m, p, q \in \mathbb{N}$, with $q \geq p$. Then we re-index $j$ and $k$ to compute a \textit{residue} for each value of $r$ such that $j$ and $k$ are as follows
$$
\begin{aligned}
& j=t m+s, \quad t=0, \ldots, p-1, \quad s=0, \ldots, m-1, \\
& k=q \ell+r, \quad \ell=0, \ldots, m-1, \quad r=0, \ldots, q-1.
\end{aligned}
$$Then the FFT (forward transform) and backward transform for $\boldsymbol{f}$ is then given by
\begin{equation}
\label{1}
\text{\textbf{Forward: }} F_{q \ell+r}=\sum_{s=0}^{m-1} \zeta_m^{\ell s} \zeta_{q m}^{r s} \sum_{t=0}^{p-1} \zeta_q^{r t} f_{t m+s},
\end{equation}
\begin{equation}
\text{\textbf{ Backward: }} f_{t m+s}=\frac{1}{q m} \sum_{r=0}^{q-1} \zeta_q^{-t r} \zeta_{q m}^{-s r} \sum_{\ell=0}^{m-1} \zeta_m^{-s \ell} F_{q \ell+r}.
\end{equation}
Computing the DFT of $\boldsymbol{f}$ then amounts to preprocessing $\boldsymbol{f}$ for each value of $r$, and computing $q$ DFTs of size $m$; no explicit zero padding is needed. The inverse transform requires $q$ FFTs of size $m$, followed by post-processing.

\subsection{Convolutions}
In mathematics and signal processing, a convolution is a mathematical operation that combines two functions to produce a third function that expresses how one function modifies the other. Specifically, it is a function that measures how much of one function is overlapping with another function as it is shifted across it. We now consider the convolution of two sequences whose inputs are $\left\{f_k\right\}_{k=0}^{M-1}$ and $\left\{g_k\right\}_{k=0}^{M-1}$ and are of finite length $M$, yielding a linear convolution with components 
\begin{equation}
    \label{4}
    (f * g)_{k}  = \sum_{p=0}^k f_p g_{k-p}, \; \; k \in \{0, \ldots, M-1\}.
\end{equation} Computing such a convolution directly requires $\mathcal{O}\left(M^2\right)$ operations, and round-off error is a problem for large $M$. This is why the convolution theorem is used instead, as it is useful in terms of complexity because it allows us to perform convolutions in the frequency domain, which can be more efficient than performing convolutions in the time domain.

\subsection{Convolution Theorem}
The convolution theorem states that the Fourier transform of a convolution of two functions is equal to the point-wise product of the Fourier transforms of the two functions. Therefore, it is preferable to use the convolution theorem, harnessing the power of the FFT to map the convolution to component-wise multiplication, and for inputs $f_k$ and $g_k$, we use the forward transform defined in \ref{dft} to obtain
\begin{equation}
\label{conv}
\sum_{j=0}^{M-1} f_j g_j \zeta_M^{-j k} = \sum_{p=0}^{M-1} \sum_{q=0}^{M-1} F_p G_q \sum_{j=0}^{M-1} \zeta_M^{(-k+p+q) j} = M \sum_s \sum_{p=0}^{M-1} F_p G_{k-p+s M},
\end{equation}
where one uses the orthogonality of the transform pair to obtain the equation. This reduces the computational complexity of convolution for \textbf{explicit padding} to $\mathcal{O}(M \log M)$ while improving numerical accuracy.

\justify
Note that if we did not explicitly pad the sequences up to length, $M$, then the output would be contaminated with aliases. This form is often used to efficiently implement numerical convolutions in a computer. It also implies that any other linear transform that has the same property of turning convolution into pointwise multiplication is essentially equivalent to the DFT up to a permutation of coefficients. This means that the DFT is unique in this sense and provides an optimal way to compute convolution efficiently. 

\begin{figure}[!ht]
    \centering
    \includegraphics[width=13cm]{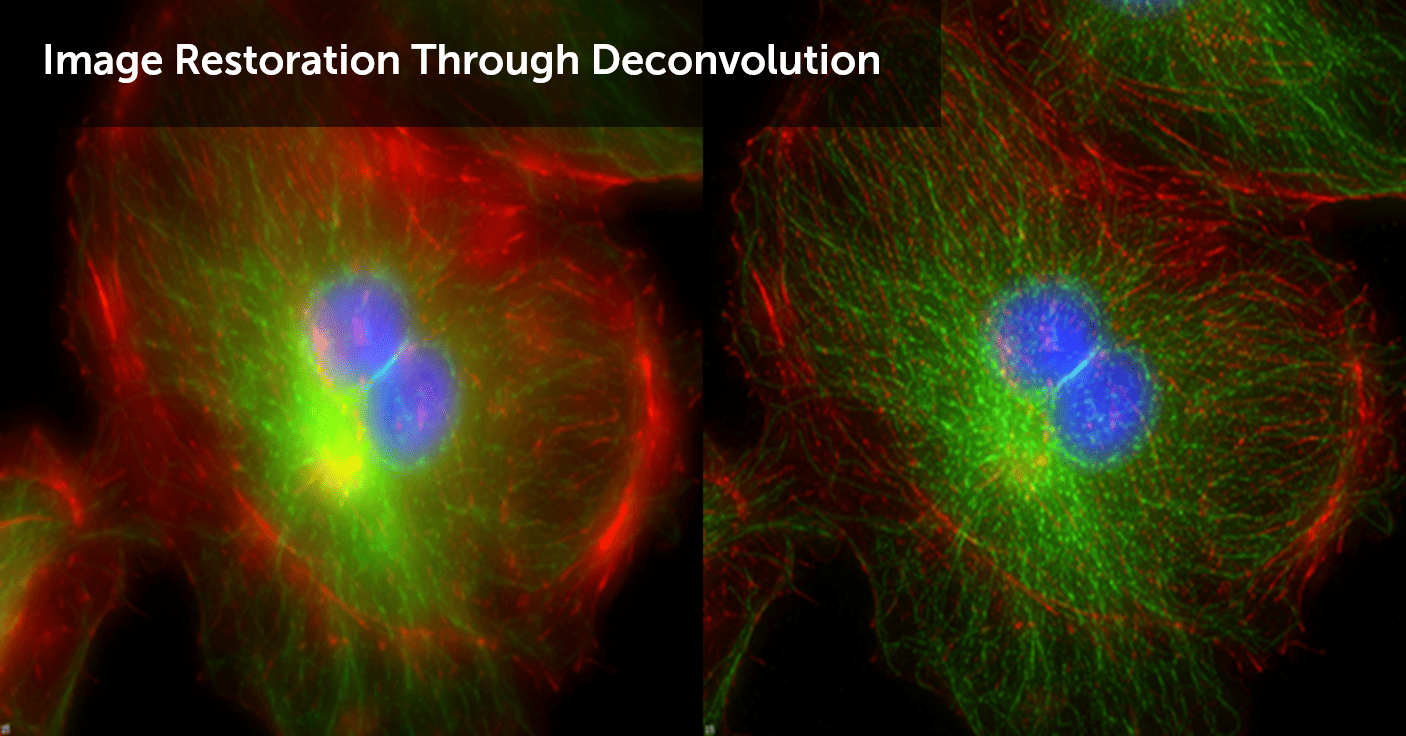}
    \caption{Initial and deconvolved images of fixed cells.}
    \label{fig:my_label}
\end{figure}

\subsection{Tuning Parameters}
\label{Parameters}
\begin{figure}
    \centering
    \includegraphics[width = 12cm]{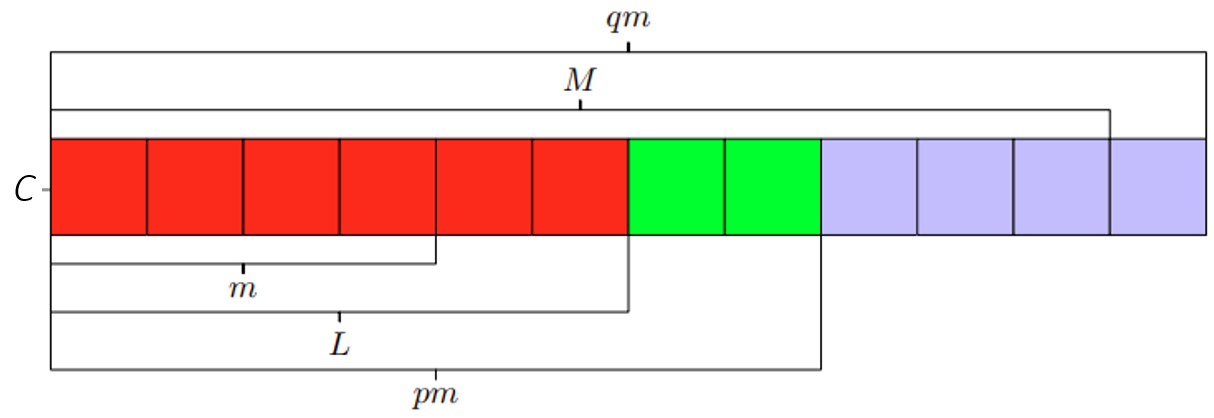}
    \caption{Tuning parameters for a 1D convolution with $L = 6$, $M = 11$, and $m = 4$, so that $p = 2$ and $q = 3$. We explicitly pad from $L = 6$ to $pm = 8$ and implicitly pad up to $qm = 12$.}
    \label{fig:1d_Parameters}
\end{figure}

We now introduce different parameters for a specific convolution for a 1D problem. Given some $m \in \mathbb{N}$, we define
$
p \doteq\left\lceil\frac{L}{m}\right\rceil, \; q \doteq\left\lceil\frac{M}{m}\right\rceil .
$
These are the smallest positive integers such that $p m \geq L$ and $q m \geq M$. So, in short, for computing the convolution between two sequences $f$ and $g$, we have 3 data-dependent parameters, which include ${\textit{L}}$ - the length of the array, ${\textit{M}}$ - the length to which we want to pad our array for the convolution, and we introduce another parameter ${\textit{C}}$ which stands for the number of copies of the FFT the transform to be computed simultaneously. If $C = 1$, we do our FFTs for each array sequentially, and if $C>1$, we can do multiple FFTs simultaneously. This is an important parameter, mostly in 2 dimensions and above, because we first compute at least $C$ FFTs in the $x$ direction and then do 1 FFT at a time in the $y$ direction before doing the convolution. Figure \ref{fig:1d_Parameters} illustrates these parameters through a 1D example. 

\justify
Now apart from the data parameters, we have 3 other important parameters, which include ${\textit{m}}$ which stands for the size of the FFTs, ${\textit{D}}$ which stands for the number of residues we compute at the time and ${\textit{I}}$ which stands for doing the FFTs in place or outplace \cite{2}. In place refers to the fact that we reuse the same input array for doing the FFTs in and putting the residues in the correct location, while outplace makes a copy of the original array and does the operations on that \textit{copied} array. 

\justify
Lastly, we can easily generalize these parameters to $n$-dimensions. To do so, we have to define these 6 parameters for each dimension, so each of these parameters, in theory, are now vectors in $\mathbb{N}^{d}$.
\begin{itemize}
\item ${\textit{L}} = (L_{1}, \ldots, L_d)$: the length of the array.
\item ${\textit{M}} = (M_{1}, \ldots, M_d)$: the length to which we want to pad our array for the convolution.
\item ${\textit{C}} = (C_{1}, \ldots, C_d)$: the number of copies of the FFT to be computed simultaneously.
\item ${\textit{m}} = (m_{1}, \ldots, m_d)$: the size of the FFTs to be used.
\item ${\textit{D}} = (D_{1}, \ldots, D_d)$: the number of residues to be computed at a time.
\item ${\textit{I}} = (I_{1}, \ldots, I_d)$: whether to do the FFTs in place or outplace.
\end{itemize}

\section{Dealiasing}
\begin{figure}[!ht]
    \centering
    \includegraphics[width = 14cm]{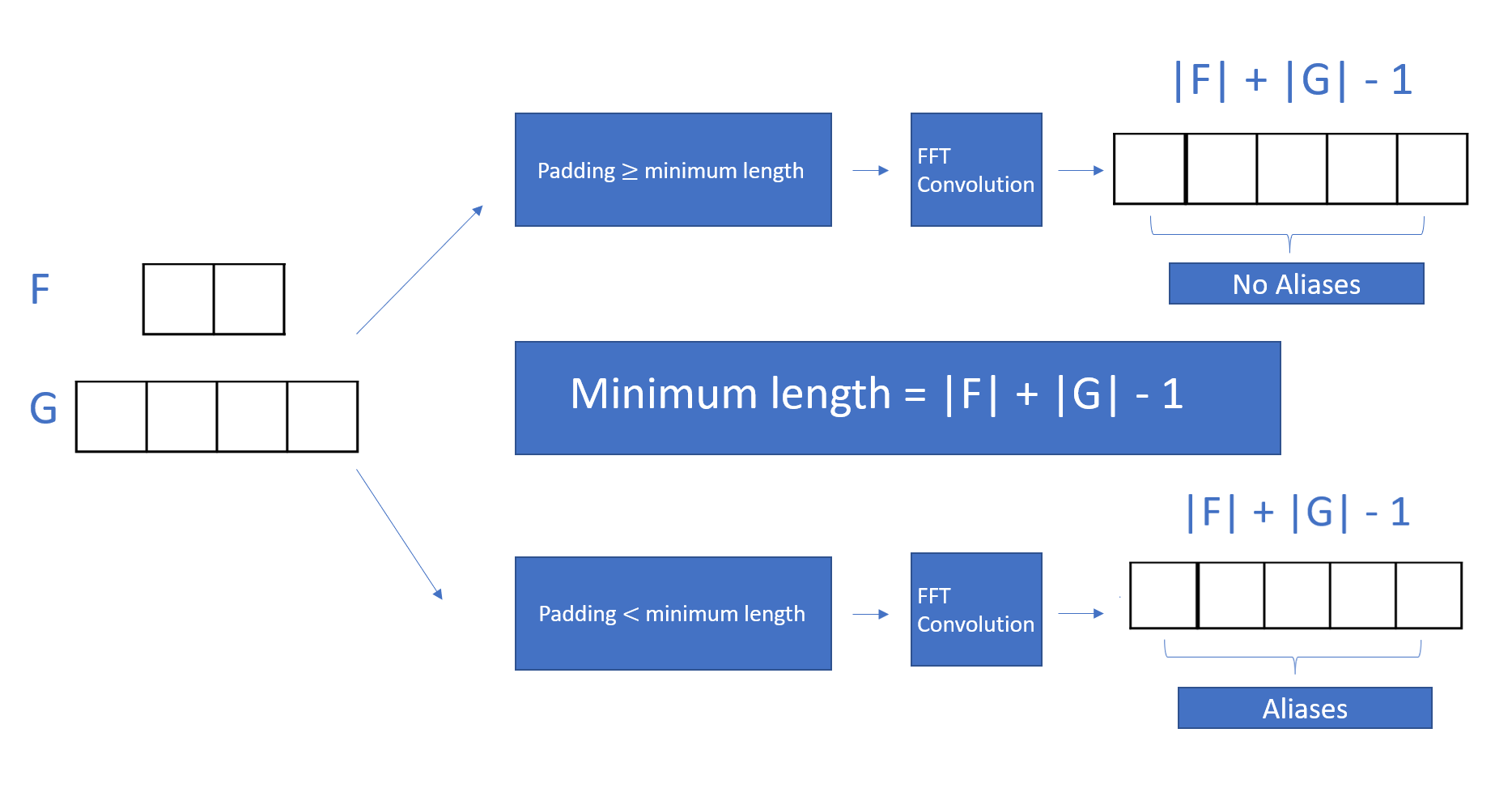}
    \caption{Linear convolutions require dealiasing of both the input arrays to a minimum size.}
    \label{fig:2}
\end{figure}

As shown in equation \ref{conv} where we obtain $M \sum_s \sum_{p=0}^{M-1} F_p G_{k-p+s M}$. The terms indexed by $s \neq 0$ are aliases, and we need to remove them to get a full linear convolution. Figure \ref{fig:2} explains why dealiasing is important for us, and in the following subsections, we will expand on the various dealiasing techniques used to solve this problem and introduce the new Hybrid Dealiasing approach.

\subsection{Explicit Dealiasing}
The classical method of dealiasing, when we want to obtain a linear convolution, is to explicitly pad zeros to both the arrays, as described in equation \ref{4}. If only the first $L$ entries of the input vectors are nonzero, aliases can be avoided by zero padding input to length $M \geq 2L - 1$. More generally, to use explicit padding, one finds a size $L\geq L_{\min}$ that the FFT is optimized for. A common choice is the next power of two, but modern FFT libraries are optimized for far more sizes, such as multiples of the form $2^a 3^b 5^c 7^d 11^e$ where $a, b, c, d, e \in \mathbb{N}$. To keep things general, we find a set of, $\tilde{L}:= \{L_1, \ldots L_d:\ L_{\min} \leq  L_j\in \mathbb{N}\ \forall\ j = 1, \ldots, d\}$ which contains values that the FFT is optimized for. The problem now can be written as follows:
$$M = \argmin_{L\in \tilde{L}} \operatorname{conv}(L).$$ where
$\operatorname{conv}(L)$ is the convolution problem, and we want to find a reasonable $M$ size for which we can efficiently use values that the FFT is optimized for.
This approach is reliable in finding an optimized size of FFT. However, while this approach guarantees an optimized size for the FFT, it still involves computing many unnecessary zero terms. This limitation motivates the development of more efficient algorithms, such as those based on fast multipole methods or low-rank approximations. Figure \ref{fig: example} shows an example in 1D.

\begin{figure}[!ht]%
    \centering
    \subfloat[\centering Explicit Padding.]{{\includegraphics[width=7.4cm]{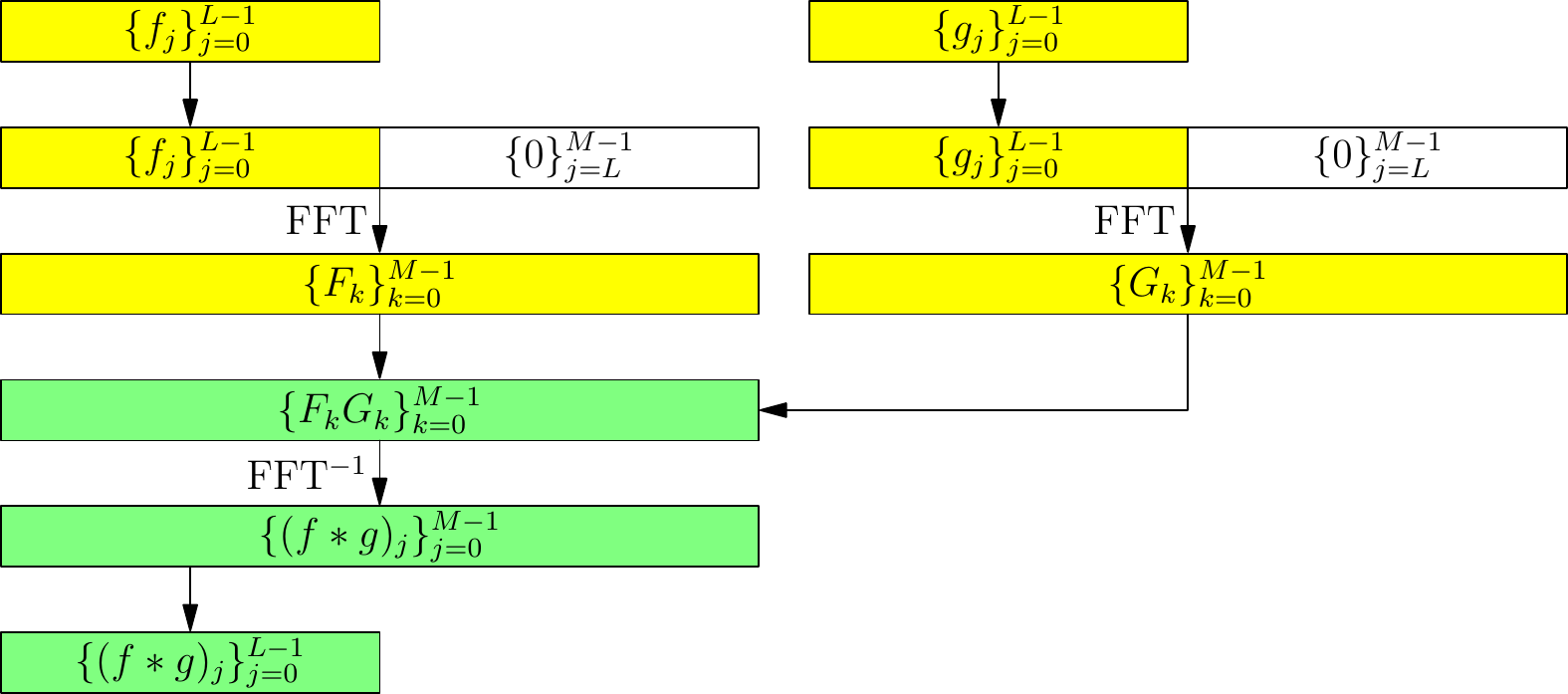} }}%
    \qquad
    \subfloat[\centering Implicit Padding.]{{\includegraphics[width=8.6cm]{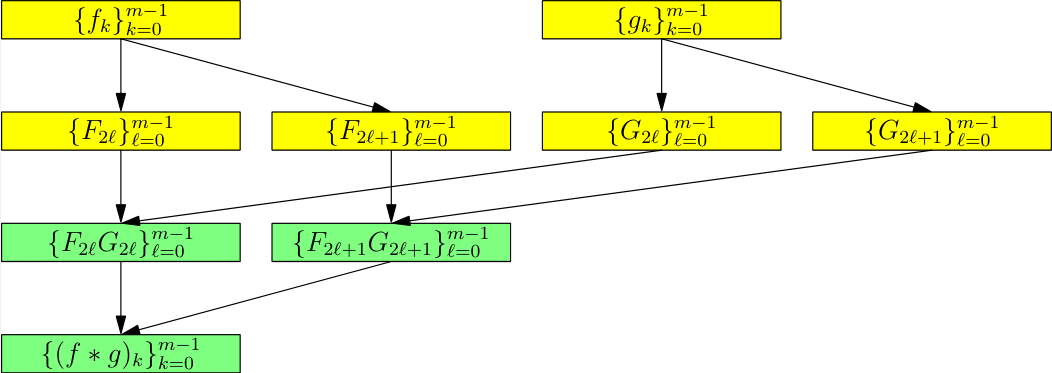} }}%
    \caption{Convolution comparison between the standard padding approaches}%
    \label{fig:example}%
\end{figure}

\subsection{Implicit Dealiasing}
An alternative approach to dealiasing with linear convolutions involves implicitly padding the input arrays rather than explicitly adding zeros as in the classical method. This approach is based on the fact that if the original array has size $L = pm$, and we want to pad it to size $M = qm$, then we can write the FFT of the padded array as $q$ FFTs of size $m$. Therefore, there is no need for explicit padding, and no computational resources are wasted on computing zero terms. The choice of $M$ depends on the possible factorization of $L$. The best value for $m$ is given by:
$$m = \argmin_{\substack{m \in \mathbb{N}\\\left\lceil\frac{L_{\min}}{m}\right\rceil m  \geq L_{\min}\\ M/m \in \mathbb{N}}} \left\lceil\frac{L{\min}}{m}\right\rceil \operatorname{conv}(m),$$
where $\operatorname{conv}(m)$ denotes the computational cost of performing a single FFT of size $m$. Once we have determined the optimal value of $m$, we can set $L_{pad} = \left\lceil\frac{M}{m}\right\rceil m$. This approach can be much faster than explicit padding, as it avoids the need to compute unnecessary zero terms. However, implicit padding is only useful if $M$ is divisible by an optimized size $m$. If $M$ is not divisible by $m$, then some zero terms still need to be computed. 
Furthermore, implicit padding relies on the size of the initial unpadded array. If $M$ is prime, it may be impossible to find a suitable value of $m$, making implicit padding infeasible. In practice, the classical method and implicit padding can be effective, depending on the problem and the specific values of $M$ and $L$.

\subsection{Hybrid Dealiasing}
Hybrid padding combines explicit padding and implicit padding in succession. We go about this as follows. Consider the minimum size we can pad our array $L_{min}$. Choose size $L_{pad}\geq L_{min}$ that can be factored into a pair $q, m \in \mathbb{N}$ with a gives a small value of $qm$. Now we explicitly pad our original array of size $L$ to size $L_{explicit} := \left\lceil\frac{L}{m} \right\rceil m$. Then we implicitly pad our array from size $L_{explicit}$ to size $L_{pad}$. This method of padding aims to use the speed of implicit padding without being constrained by the factors of the size of the original array. The main problem with hybrid padding is how one efficiently determines $L_{pad}$. 

\subsection{Powers Of 2 Example}

For a more general case, we choose $M = 2L_2 = 2^\alpha L_1$ for some $\alpha > 0$. In case we have odd-length arrays, pad the arrays explicitly so it fits the above constraint. The convolution is given by (forward transform)

\[F_j = \sum_{j = 0}^{M-1} \zeta_{M}^{kj} f_j \text{ and } G_j = \sum_{j = 0}^{M-1} \zeta_{M}^{kj} g_j,\]

\justify
and $f_j = 0  \text{ for } k \geq L_1$ and $g_j = 0  \text{ for } k \geq L_2$. Now instead of dividing the sequence into 2 arrays, we divide them into $2^\alpha$ arrays ($\frac{2L_2}{L_1} = 2^\alpha$) for $f$ and 2 arrays as usual for $g$. This is given by
\[F_{2^\alpha k} = \sum_{j = 0}^{L_1-1} \zeta_{L_1}^{kj} f_j, \; F_{2^\alpha k+1} = \sum_{j = 0}^{L_1-1} \zeta_{L_1}^{kj} \zeta_{M}^{j} f_{j},  \;F_{2^\alpha k+2} = \sum_{j = 0}^{L_1-1} \zeta_{L_1}^{kj} \zeta_{M}^{2j} f_{j},\; \hdots F_{2^\alpha k + 2^{\alpha-1}} = \sum_{j = 0}^{L_1-1} \zeta_{L_1}^{kj} \zeta_{M}^{2^{\alpha-1}j} f_{j},\]
and 
\[G_{2k} = \sum_{j = 0}^{L_2-1} \zeta_{L_2}^{kj} g_j \text{ and } G_{2k+1} = \sum_{j = 0}^{L_2-1} \zeta_{L_2}^{kj} \zeta_{M}^{j} g_{j}.\]
\justify
We can do this by multiplying the above arrays in the following manner;

\[A_0 = \{F_0[0],F_2[0], \cdots, F_{2^{\alpha}-2}[0], F_0[1],F_2[1], \cdots F_{2^{\alpha}-2}[1], \cdots, F_{0}[L_1 - 1], \cdots F_{2^{\alpha}-2}[L_1 - 1]\},\]
\[A_1 = \{F_1[0],F_3[0], \cdots,F_{2^{\alpha}-1}[0], F_1[1],F_3[1], \cdots,F_{2^{\alpha}-1}[1], \cdots, F_{1}[L_1 - 1], \cdots F_{2^{\alpha}-1}[L_1 - 1]\}.\]
\justify
Now use these two arrays ($A_0 \text{ and } A_1$) to multiply the following arrays $G_0$ and $G_1$ respectively;

\[H_0 = A_0 \odot G_0 \rightarrow \operatorname{FFT^{-1}} \text{ and } H_1 = A_1 \odot G_1 \rightarrow \operatorname{FFT^{-1}}. \]
and now, taking the inverse Fourier transform of the above arrays, we get each array $H_0$ and $H_1$ are of length $L_2$, respectively. The orthogonality of this transform pair follows from
$$
\sum_{j=0}^{M-1} \zeta_{M}^{\ell j}= \begin{cases}M & \text { if } \ell=s M \text { for } s \in \mathbb{Z} \\ \frac{1-\zeta_{M}^{\ell M}}{1-\zeta_{M}^{\ell}}=0 & \text { otherwise. }\end{cases}
$$
\justify
Therefore then, by the convolution theorem, we get
$$\begin{aligned}
	\sum_{j=0}^{M-1} F_{j} G_{j} \zeta_{M}^{-j k} &= \sum_{j=0}^{M-1} \zeta_{M}^{-j k}\left(\sum_{p=0}^{M-1} \zeta_{M}^{j p} f_{p}\right)\left(\sum_{q=0}^{M-1} \zeta_{M}^{j q} g_{q}\right) = \sum_{p=0}^{M-1} \sum_{q=0}^{M-1} f_{p} g_{q} \sum_{j=0}^{M-1} \zeta_{M}^{(-k+p+q) j} =M \sum_{s} \sum_{p=0}^{M-1} f_{p} g_{k-p+s M} .
\end{aligned}$$

\justify
Lastly, the backward transform of the convolution is given by

\[Mf_k = \sum_{j = 0}^{L_2-1} \zeta_{L_2}^{-kj}H_{2j} + \zeta_{M}^{-k} \sum_{j = 0}^{L_2-1} \zeta_{L_2}^{-kj}H_{2j+1} \; \forall k \in \{0,1, \cdots L_2 - 1\}.\]

\subsection{Final Algorithm For One Dimension}
In order to compute the convolution of two arrays $f$ and $g$ of different sizes $L_f$ and $L_g$ that both need to be padded to a size $M$, we can use. However, this requires some assumptions to be made. We assume that $L_f \leq L_g$ (this assumption is made for convenience, as we can always call the larger array $g$) and uses FFTs of size $m$. Since the arrays are of different sizes, we can use different values of $m$. We assume that $f$ uses a value of $m$ that is smaller than or equal to the $m$ used for $g$. The values $L_f, L_g$, and $M$ are all problem-dependent.
To carry out , we introduce two integers $m, \Lambda \in \mathbb{N}$ and define the following values:
$$
\begin{aligned}
p_f :=\left\lceil\frac{L_1}{m}\right\rceil, 
p_g :=\left\lceil\frac{L_2}{\Lambda m}\right\rceil, 
q_g :=\left\lceil\frac{M}{\Lambda m}\right\rceil,
q_f :=\frac{q_g \Lambda m}{m}=\Lambda q_g .
\end{aligned}
$$
Note that because $\Lambda m \geq m$, we have $q_g \leq q_f$.
For given integers $r_g \in\left[0, q_g-1\right]$ and $\ell_g \in[0, \Lambda m-1]$, define,
$$
\begin{aligned}
& \ell_f=\left\lfloor\frac{\ell_g}{\Lambda}\right\rfloor, \lambda=\ell_g \bmod \Lambda, r_f=q_g \lambda+r_g.
\end{aligned}
$$
Then $$
\begin{aligned}
q_g \ell_g+r_g & =q_g\left(\Lambda \ell_f+\lambda\right)+r_g =\Lambda q_g \ell_f+\lambda q_g+r_g =q_f \ell_f+r_f .
\end{aligned}
$$
This means that for each residue $r_g$ of $g$, we can compute the corresponding residues $r_f=q_g \lambda+r_g$ of $f$ for each $\lambda=0, \ldots, \Lambda-1$. We can then perform all necessary multiplication and take the inverse transform. We can reuse the memory to compute the next residues. After one has computed $H$, one can take the inverse transform and accumulate appropriately. Using this, we can efficiently compute the convolution of two arrays of different sizes without explicitly zero-padding them. This approach is useful when the sizes of the arrays are not known in advance or when there are computational constraints. The pseudo-codes for the forward and backward routines are shown in \ref{algo_forward} and \ref{algo_backward}, and the main convolution in \ref{mainconvolve}. \cite{2} shows the Hybrid Dealiasing algorithm for 1-dimensional convolution.

\subsection{Final Algorithm for $N$ Dimensions}

\begin{figure}[!ht]
    \centering
    \includegraphics[width = 18cm]{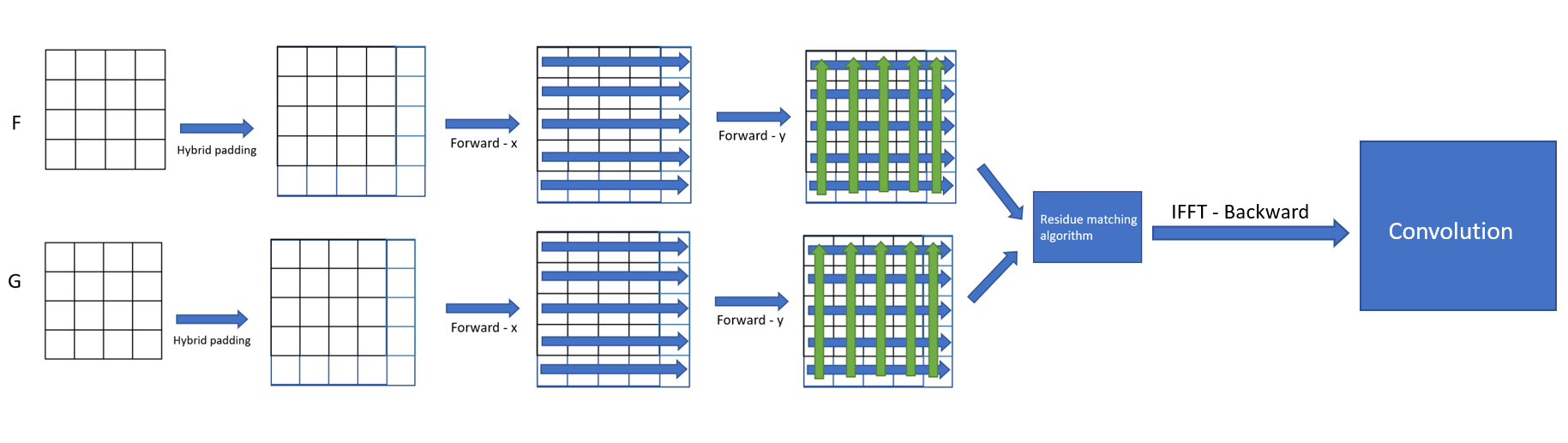}
    \caption{Hybrid Dealiasing algorithm for 2D.}
    \label{fig:3}
\end{figure}

To generalize the  method to $N$ dimensions, we can introduce the same parameters for each dimension. Let $f$ and $g$ be two $N$-dimensional arrays of sizes $L_{f,1}, \ldots, L_{f,N}$ and $L_{g,1}, \ldots, L_{g,N}$, respectively, and let $M_{1}, \ldots, M_{N}$ be the size of the padded arrays. Let $m_{1}, \ldots, m_{N}$ and $\Lambda_{1}, \ldots, \Lambda_{N}$ be integers that define the sizes of the FFTs to be used in each dimension. We assume that $L_{f,i} \leq L_{g,i}$ for all $i \in {1, \ldots, N}$, and that $m_{i}$ is smaller than or equal to $\Lambda_{i}$ for all $i \in {1, \ldots, N}$. For each dimension, $i$, we define the following values:
$$
\begin{aligned}
p_{f, i} :=\left[\frac{L_{f, i}}{m_i}\right],
p_{g, i} :=\left[\frac{L_{g, i}}{\Lambda_i m_i}\right], 
q_{g, i} :=\left[\frac{M_i}{\Lambda_i m_i}\right],
q_{f, i} :=\frac{q_{g, i} \Lambda_i m_i}{m_i}=\Lambda_i q_{g, i} .
\end{aligned}
$$
Note that because $\Lambda_{i} m_{i} \geq m_{i}$, we have $q_{g,i} \leq q_{f,i}$. For each dimension $i$, let $r_{g,i} \in\left[0, q_{g,i}-1\right]$ and $\ell_{g,i} \in[0, \Lambda_{i} m_{i}-1]$. We can then define
$$
\begin{aligned}
& \ell_{f, i}=\left\lfloor\frac{\ell_{g, i}}{\Lambda_i}\right\rfloor, \lambda_i=\ell_{g, i} \bmod \Lambda_i, r_{f, i}=q_{g, i} \lambda_i+r_{g, i} .
\end{aligned}
$$
Then we have
$$
\begin{aligned}
q_{g, i} \ell_{g, i}+r_{g, i} & =q_{g, i}\left(\Lambda_i \ell_{f, i}+\lambda_i\right)+r_{g, i} =\Lambda_i q_{g, i} \ell_{f, i}+\lambda_i q_{g, i}+r_{g, i} =q_{f, i} \ell_{f, i}+r_{f, i} .
\end{aligned}
$$
This means that for each residue $r_{g,i}$ of $g$ in dimension $i$, we can compute the corresponding residues $r_{f,i}=q_{g,i} \lambda_{i}+r_{g,i}$ of $f$ for each $\lambda_{i}=0, \ldots, \Lambda_{i}-1$ in that dimension. Then we call the forward transform (FFT) $N – 1$ dimension routine recursively along the $N$th axis. We
then multiply the two arrays in Fourier space using the Residue matching algorithm and call the backward transform $N – 1$ dimension routing recursively
along the $N$th axis. Figure \ref{fig:3} shows an example of the 2D convolution algorithm.

\section{Optimization Search}
Finding the best parameters for numerical algorithms is essential because it can significantly impact the performance and accuracy of the algorithm. This is because numerical algorithms rely on input parameters to make decisions and carry out computations, and choosing these parameters can significantly affect the results. For example, consider an optimization algorithm that tries to minimize a function. The choice of optimization parameters, such as the learning rate, convergence criteria, and regularization strength, can significantly affect the convergence rate, final solution quality, and generalization ability of the algorithm. Similarly, in machine learning algorithms such as neural networks, the choice of hyperparameters such as learning rate, batch size, and the number of layers can significantly affect the training process and the model's performance on unseen data.
\justify
In summary, finding the best parameters for numerical algorithms is crucial for achieving the best possible performance and accuracy. Furthermore, it can differentiate between a successful and unsuccessful algorithm or model. We discuss two basic approaches we developed to find the best hyperparameters for our algorithms, as described in \ref{Parameters}. To emphasize, the specific parameters we want to search over are $L$, $C$, $m$, $D$ and $I$. 

\subsection{Problem}
Our problem is simple: Given a large $N$ dimensional input, we want to estimate the best convolution parameters. For simplicity, we restrict ourselves to 2 dimensions and look for rectangles whose parameters can be used to estimate the square's parameters. More formally put, consider a square of size $L_{xy}$ where $x = y$ for squares (for simplicity, we denote this by $L_x$) and $x \neq y$ for rectangles. We also have another parameter, i.e., the size to which we pad which is denoted by $M_{xy}$ and for a square, just denoted by $M_x$. We fix a square of size $L_{x}$ and $M_{x} = 2*L_{x}$ and find rectangles with the same width $x$ fixed but varying the height, $y$ where the parameters of such a rectangle $L_{xy}$ and $M_{xy} = 2*L_{xy}$ can be used to estimate the square's parameters. The goal is to find the \textit{skinniest} rectangles to make this possible.

\subsection{Naive Approach – {A0}}
${A0}$ used a grid search to find these optimal parameters, which is inefficient as it scales with the number of parameters. We have 5 parameters for each dimension $d_i, \; i \in \{1,2 \ldots N\}$. If we want to find the optimal parameters, we would have to go through $N^5$ combinations, and hence our algorithm would take $\mathcal{O}(N^5)$. In general, if we have $K$ parameters, our algorithm would grow asymptotically $\mathcal{O}(N^K)$ if we did them in a nested manner. This is a popular method called \textit{grid search}. Grid search involves defining a range of values for each hyperparameter and then exhaustively searching through all possible combinations of hyperparameter values. This process generates a grid of hyperparameter combinations, and each combination is then used to train and evaluate an algorithm/model. The performance of each model/algorithm is then recorded and compared, and the hyperparameter combination that yields the best performance is selected as the optimal set of hyperparameters. However, for our problem, we can exploit a trick that reduces the computational complexity significantly from $\mathcal{O}(N^5)$ to $\mathcal{O}(N)$ as we notice that we can search for the optimal hyperparameters independently for each dimension due to the very nature of how a multidimensional FFT is calculated and this is optimal for any dimension. 

\justify
We can search for the optimal hyperparameters independently for each dimension of the multidimensional FFT based on the FFT's separability property. The separability property of the FFT means that we can perform the FFT on each dimension independently and then combine the results to obtain the final multidimensional FFT. This property arises because the multidimensional FFT can be expressed as a sequence of one-dimensional FFTs. As a result, the optimal hyperparameters for each dimension can be found independently without affecting the performance of the overall multidimensional FFT. This is because the optimal hyperparameters for each dimension only affect the computation of the one-dimensional FFT for that dimension and do not affect the computation of the other dimensions. Therefore, by exploiting the separability property of the FFT, we can significantly reduce the computational complexity of searching for the optimal hyperparameters of a multidimensional FFT and achieve optimal performance for any number of dimensions.

\subsection{Experience – {A1}}
In order to create a more optimized and efficient hyperparameter tuning algorithm, we have developed ${A1}$, which is designed to be better than ${A0}$ and is capable of searching over a broader set of parameters to determine the fastest algorithm for a given problem empirically. To achieve this, we consider a set of forward routines $\textbf{F}$, which are algorithms specialized for some instances depending on parameters or the structure of the data and the convolution required. We also have a set of backward routines $\textbf{B}$, which we want to optimize over. However, since the backward and forward transforms are related, we can optimize over the forward routines only and use the same optimal parameters for the backward routine. It is important to emphasize that our benchmark is the timing of the forward routine ((Let us denote this score by ${\textit{R}}$). In order to make $\textbf{A1}$ computationally efficient, generalizable, and scalable for future routines, we have developed several approaches:
\begin{enumerate}
    \item Firstly, we have adopted a more efficient approach called \textit{random search}, which is heavily used in machine learning. This approach randomly samples a smaller set of values (denoted by $N_1$) from each parameter set (according to a distribution) where $N_1 < N$. This approach is more efficient than grid search, as it explores only a small subset of the parameter space but still finds the optimal hyperparameters.

\item Secondly, we have explored the use of a tuning algorithm called \textit{Bayesian search}, which is gaining popularity in ML. This approach optimizes its parameter selection in each round according to the previous round $\textbf{\textit{R}}$. This means that if certain combinations of parameters give good results in ${\textit{R}}$, then the algorithm will continue to explore that space instead of exploring a space which gives poor results. Although this approach has a slower learning rate than random search, it is more likely to find the ideal set of parameters.

\item Lastly, we have developed efficient heuristics by understanding the various relationships between parameters and avoiding recomputing the parameters for specific data. For example, if the optimal set of parameters for a rectangle and a square are the same, then we can compute the parameters on the rectangle and reuse them in the square. We can reduce the overall computation time required for hyperparameter tuning by leveraging such relationships between parameters.
\end{enumerate}

\justify
Our ${A1}$ algorithm is designed to be more efficient, generalizable, and scalable than ${A0}$. Using a combination of random search, machine learning and heuristics, we can efficiently search over a broader set of parameters and determine the fastest algorithm for a given problem.

\subsection{Machine Learning}
Machine learning is a rapidly growing field that has gained widespread attention in recent years due to its potential to revolutionize how we approach data analysis and prediction. At its core, machine learning involves developing algorithms and models that can automatically identify data patterns, learn from them, and make predictions or decisions based on that learning. One area where machine learning has the potential to make a significant impact is in numerical algorithms, which are used to solve complex mathematical problems that are often too difficult or time-consuming to solve by hand. By using machine learning techniques to develop more efficient and accurate numerical algorithms, we can improve our ability to solve various problems in physics, engineering, finance, and more. Specifically, machine learning can help us in numerical algorithms by enabling us to:
\begin{enumerate}
    \item Develop more accurate models for predicting outcomes based on data
\item Optimize algorithms for faster and more efficient computation
\item Identify and handle complex or non-linear relationships between variables
\item Handle large datasets and complex data structures more effectively
\end{enumerate}
Overall, machine learning can enhance our ability to solve complex mathematical problems and make more accurate predictions, ultimately helping us better understand how these parameters are correlated and use that to build complex algorithms that would efficiently let us predict optimal parameters.

\section{Methodology}
We now provide the reader with a clear understanding of how the problem was tackled; the steps are taken to reach a solution, and the rationale behind each decision. In this section, we will discuss the different methods and techniques used to optimize hyperparameters and improve the efficiency of the forward transform algorithm. We will describe the advantages and disadvantages of each approach and explain how they were implemented to achieve the desired results in the subsequent section.
\begin{figure}
    \centering
    \vspace{-0.4cm}
    \includegraphics[width = 15cm]{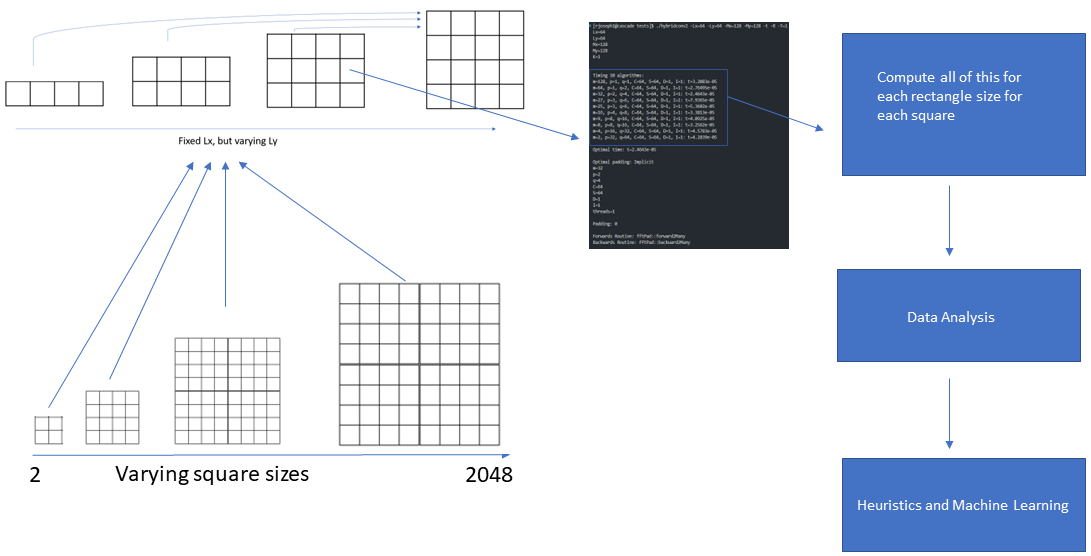}
    \label{fig:my_label}
    \caption{Complete machine learning pipeline.}
\end{figure}  

\subsection{Data Generation}
Initially, we planned to search over all rectangles $L_{xy}$ where $y \in \{2, \ldots x\}$ for each fixed square $L_x$. However, this approach proved time-consuming, and the results were random, needing discernible patterns. Subsequently, we modified our approach and restricted our search to squares and rectangles of size $2^a 3^b 5^c 7^d \text{ where } a,b,c,d \in \mathbb{N}$. This new approach yielded significantly better results than the original method. Moreover, searching over these specific sizes was more efficient and faster regarding hardware computation. These are the complete steps that we followed for caching the dataset

\begin{enumerate}
    \item For each fixed square $L_x$ in $\operatorname{range}(2,3,5,7)$ we search over all rectangles $L_{xy}$ where $y \in \operatorname{range}(2,3,5,7)$ where $\operatorname{range}(2,3,5,7)$ are rectangles/squares of size $2^a 3^b 5^c 7^d$ and $a,b,c,d \in \mathbb{N}$.
\item For each rectangle $L_{xy}$, we use an optimizer developed in \cite{2} to benchmark various timings. We repeat the optimizer runs twice and cache all the results. Specifically, we record the $x$ direction FFT $m$ and the associated time values.
\end{enumerate}

\subsection{Data Analysis: $\epsilon$ Threshold}
We now checked for any patterns by tackling our goal by introducing this threshold factor $ \epsilon $. We begin by searching for the best parameters for the square, $L_x$ denoted by $m_{s}^{*}$ and $time_{s}^{*}$. We now create another threshold factor, which we call $\gamma$ where $\gamma(\epsilon) = time_{s}^{*} + \epsilon * (time_{s}^{*} - time_{s}^{worst})$. Then again, as usual, we search over all the rectangles $L_{xy}$ where $y \in \{2, \ldots x\}$. For each rectangle, we get a sequence of $m$ values denoted by $m_{r}^{i}$ and time values $time_{r}^{i}$ where $i \in \{1, 2 \ldots 13\}$. We then check if $m_{s}^{*}$ is within the $\gamma(\epsilon)$ factor of the sequence of $m_{r}$ produced by that rectangle. We looked for the smallest rectangle size to achieve this for different $\epsilon$ values. 

\subsection{Regression vs Classification}
Lastly, we attempted to develop several machine learning algorithms to predict the rectangle $L_{xy}$ from the $L_{xy}$ value, $M_{xy}$ value, and a specific $\epsilon$ value. We treated this as a regression task since we wanted to estimate an integer and compared various standard regression algorithms such as Neural Networks, Logistic Regression, Lasso/Linear/Ridge Regression, Elastic Net regression, Random Forest Regressor, Decision tree regressors, Extra Trees, and ensembles of regression models.

\begin{enumerate}
    \item To fit these models to the dataset, we first cached all the data into a pandas data frame, which includes three features: $L_{xy}$, $M_{xy}$, and $\epsilon$, with the target output being $L_{y}$, where $L_{y}$ represents the smallest rectangle within $\epsilon$ factor of the large square. We then split the dataset by a random 0.3 split to evaluate our performance. The training set contained 70\% of the data, while the test set contained the remaining 30\%.

\item The reason for splitting the data into training and test sets is to simulate how the model would perform on unseen data. Evaluating the model's performance on a separate test set allows us to estimate how well it will likely perform on new, unseen data. It is crucial to note that the test set should not be used during training to avoid biasing the performance metrics.

\item Next, we fit each model and evaluated its performance using 5-fold cross-validation. This technique involves dividing a dataset into five subsets of equal size, using four subsets to train the model and the remaining subset to test the model's performance. This process is repeated five times, and the final performance metrics are calculated by taking the average of the metrics across all five iterations. Using 5-fold cross-validation provides a more robust estimate of model performance and helps to reduce the risk of overfitting. We choose the best model out of this and proceed to test the model on the test dataset.
\end{enumerate}

\justify
We repeated the same process using various classification algorithms to switch to a multi-class classification task. Instead of being a regression task, we treated this as a multi-class classification task. We evaluated the performance of several standard classification algorithms and selected the best model based on our evaluation. 

\section{Results and Discussions}

\subsection{1D/2D Convolutions}
\begin{figure}[!ht]%
    \centering
    \subfloat[\centering Fixed-size array $f$ of length 2.]{{\includegraphics[width=8cm]{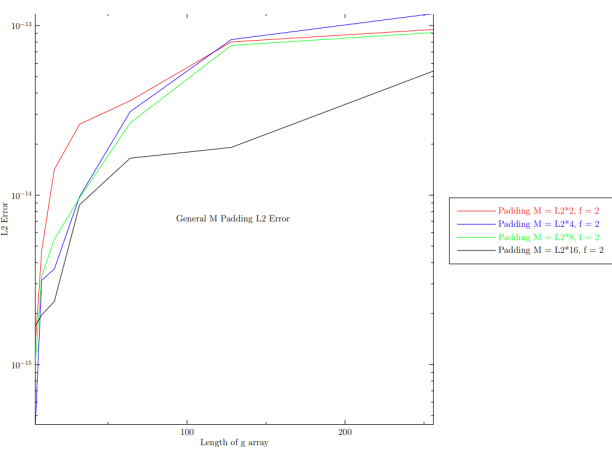} }}%
    \qquad
    \subfloat[\centering Fixed-size array $f$ of length 12.]{{\includegraphics[width=8cm]{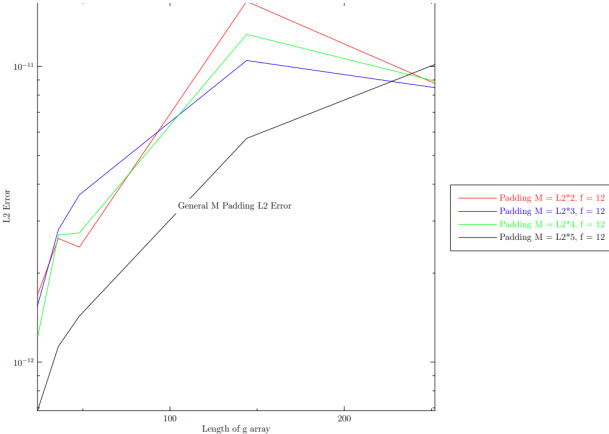} }}%
    \caption{Convolution comparison where the size of $g$ get varied against the $L_2$ error for different padding size. }%
    \label{fig:plotl2}%
\end{figure}

We offer unique solutions that expand hybrid dealiasing to issues with variable input sizes, reducing memory and computation time. We expanded this to multi-convolution (1D/2D), applying the convolution to a data sequence $N$ arrays rather than just two. We construct the first hybrid dealiasing solution quicker than typical explicit/implicit dealiasing approaches. We compare various padding errors, too, as shown in Figure \ref{fig:plotl2} by fixing the size of $f$.

\subsection{$\epsilon$ Varying Heuristics}

\begin{figure}[!ht]
    \centering
    \includegraphics[width = 18cm]{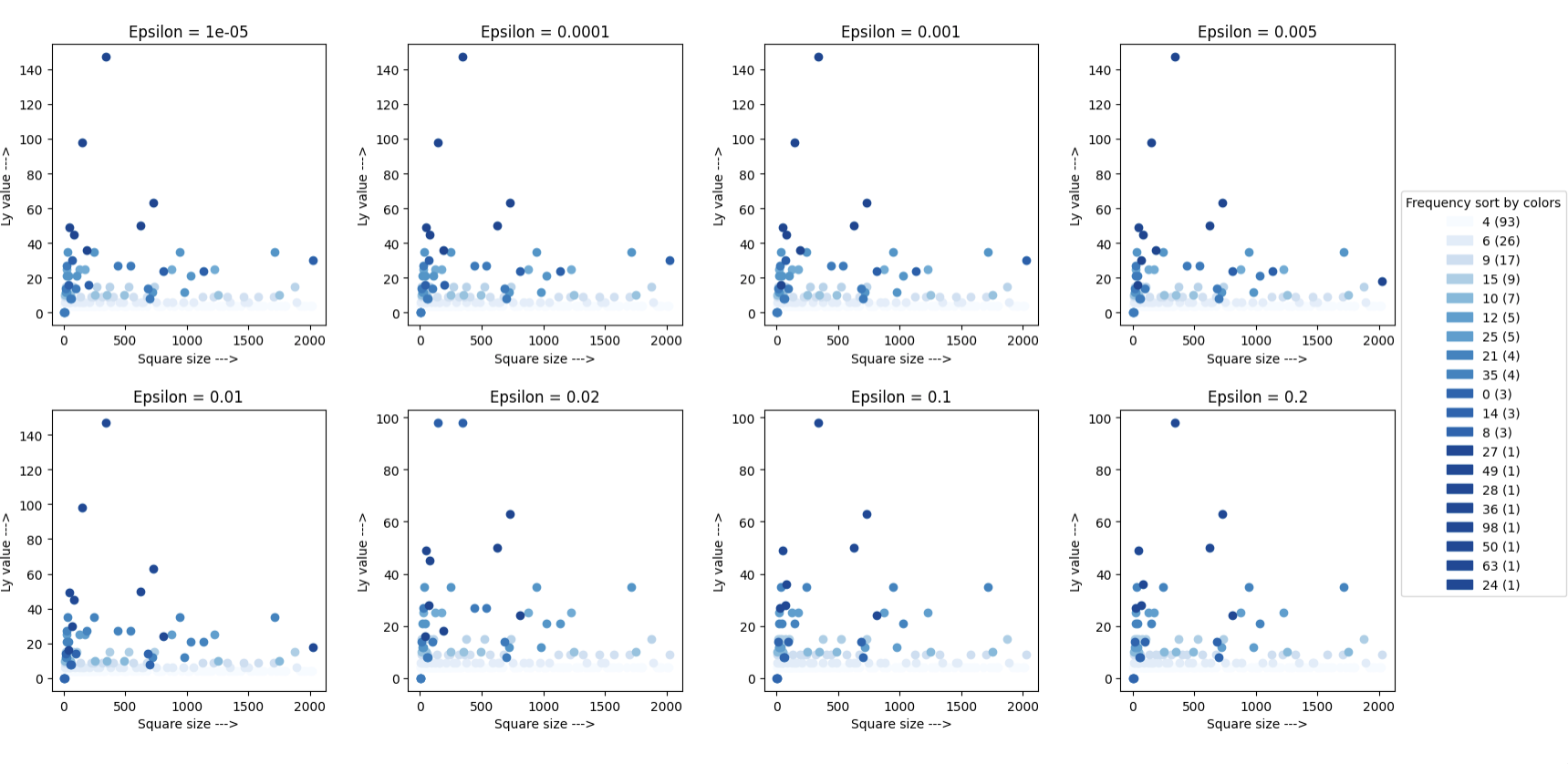}
    \caption{Various $\epsilon$ factor plots. On the x-axis, we vary the square size while the y-axis represents the smallest rectangle size found ($L_y$).}
    \label{fig:epsilon1}
\end{figure}    

In Section 5.4, we introduced $\epsilon$ as the threshold at which the optimal sample complexity $m_{s}^{*}$ is within a factor of $\epsilon$ of the smallest sample complexity $m_{r}$. Figure \ref{fig:epsilon1} presents some intriguing findings, the foremost being that an increase in $\epsilon$ leads to tighter bounds on $L_y$ values. This can be easily explained as an increase in $\epsilon$ results in more values to compare against, which makes it easier to identify smaller values. Notably, the upper bound of 140 suggests that even with a small value of $\epsilon$ such as $0.005\%$, we can obtain highly accurate approximations of the parameters of the large squares using significantly smaller (narrower) rectangles (maximum height of $L_y = 140$).
 
 \subsection{Prediction Of $L_y$}
   \begin{figure}[!ht]
    \centering
    \includegraphics[width = 12cm]{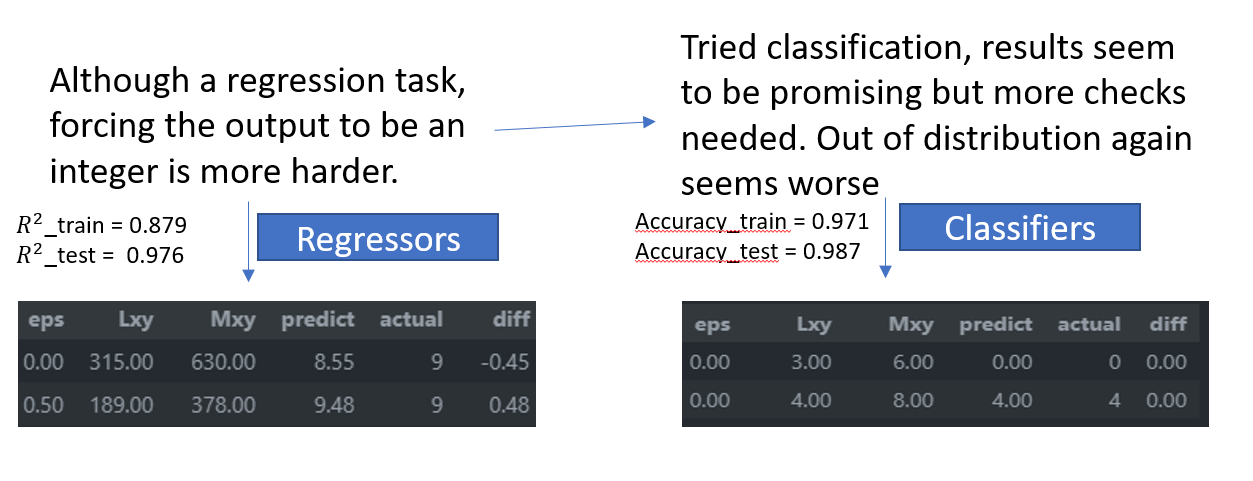}
    \caption{Comparison between the regression and classification estimators.}
    \label{fig:ml}
\end{figure} 

\justify
Our model's regression estimators could have performed better on out-of-distribution datasets. This was because we tried to estimate the size of rectangles as continuous parameters and then constrain them to integer values, resulting in significant estimation errors. In contrast, classification-based approaches better deal with discrete values and provide a more accurate approximation of the underlying data. When using regression models, the predicted rectangle size values are continuous and must be rounded to integer values. This can lead to significant errors, as even slight deviations from true integer values can significantly impact the model's overall accuracy, as shown in Figure \ref{fig:ml}.
In contrast, classification-based approaches categorize data into discrete groups, providing a more accurate approximation of the underlying data. Therefore, classification-based estimators fit the dataset better because they are better suited to the nature of the problem. Additionally, classification-based models can handle the problem of imbalanced classes more effectively than regression-based models. Figure \ref{fig:ml} illustrates the differences and the results between the two approaches.

 \subsection{Heuristic Guarantees}
We finally consider a different approach from the $\epsilon$ threshold, where instead of finding the smallest $L_{xy}$ rectangle for a given square, we instead compare the rectangle of height 32 ($L_y = 32$) and use that rectangles best $m$ values to compare it to the square's best $m$ values. More formally put, for any square $L_x$, we fix the rectangle height to $L_y = 32$, find the best $m$ values ($m_{r}^{*}$), check to see if $m_{r}^{*}$ is in the sequence of values produced by the square denoted by $m_{s}^{i}$ and time values $time_{s}^{i}$ where $i \in \{1, 2 \ldots 13\}$. If it is in that sequence, let us denote that index by $j$, we compare $time_{s}^{j}$ and $time_{s}^{*}$ as well as compare $m_{s}^{j}$ and $m_{s}^{*}$ to see how close to each other they are.

\subsubsection{$m$ Values}
In figure \ref{fig:plotsfinal}, we observe the presence of three distinct rays of $m$ values. This interesting find suggests limiting our search to only these three rays of $m$ values to approximate larger squares. The intriguing question is, what could cause these rays to appear in the first place? One plausible conjecture is that each of these rays corresponds to a different type of padding ratio, as alluded to in \cite{2}. In other words, each of these rays of $m$ values may be related to a specific padding scheme used when padding is required to adjust the input size of the convolution algorithm. Therefore, by focusing our search on these specific rays, we can optimize our approximation algorithm and significantly reduce the computational complexity required to obtain accurate results. Future research in this area could explore the underlying mechanisms that give rise to these rays and investigate how this observation can be generalized to other types of convolution problems.

\subsubsection{Time Values}
\begin{figure}[!ht]%
    \centering
    \subfloat[\centering $m$ values comparison]{{\includegraphics[width=8.3cm]{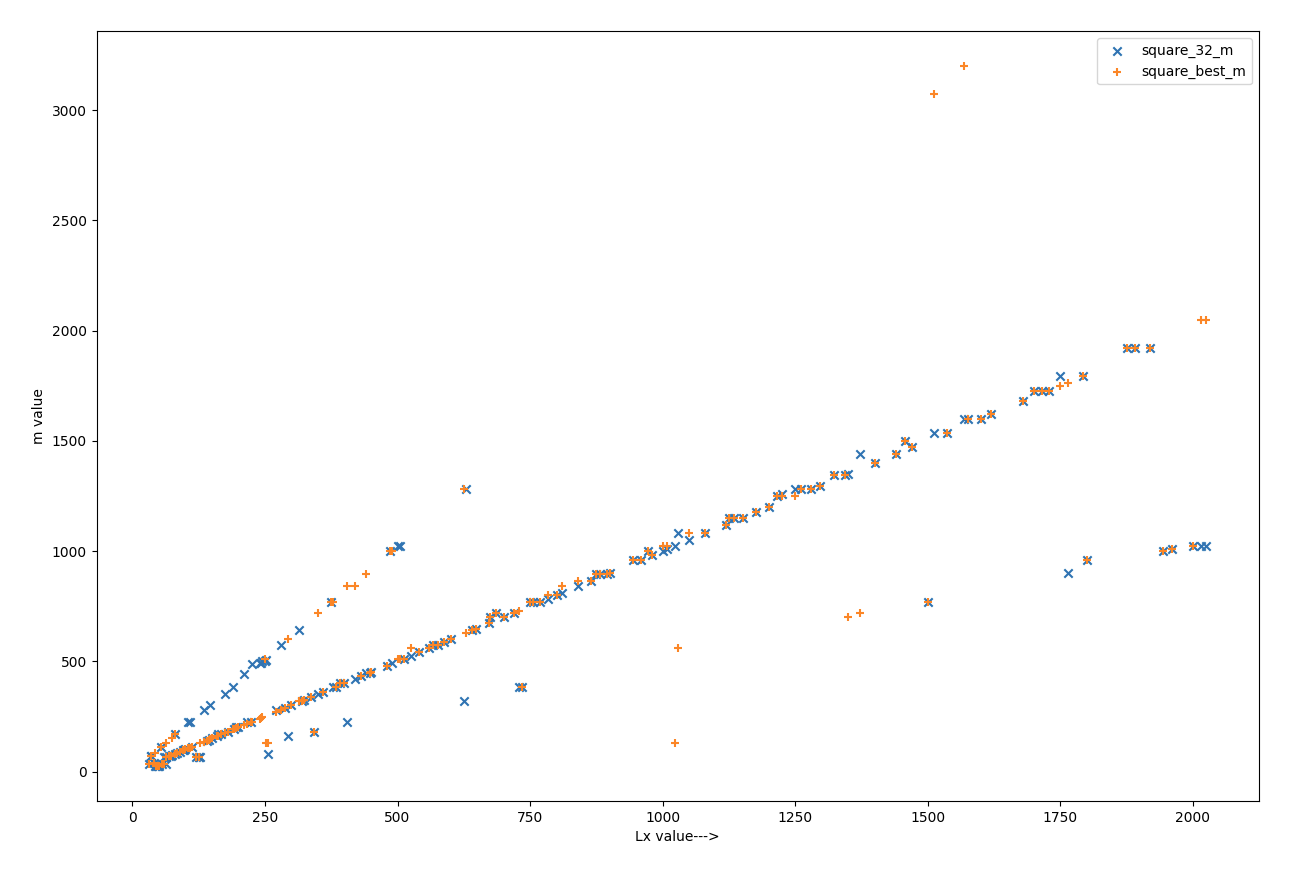} }}%
    \qquad
    \subfloat[\centering time values comparison]{{\includegraphics[width=8.3cm]{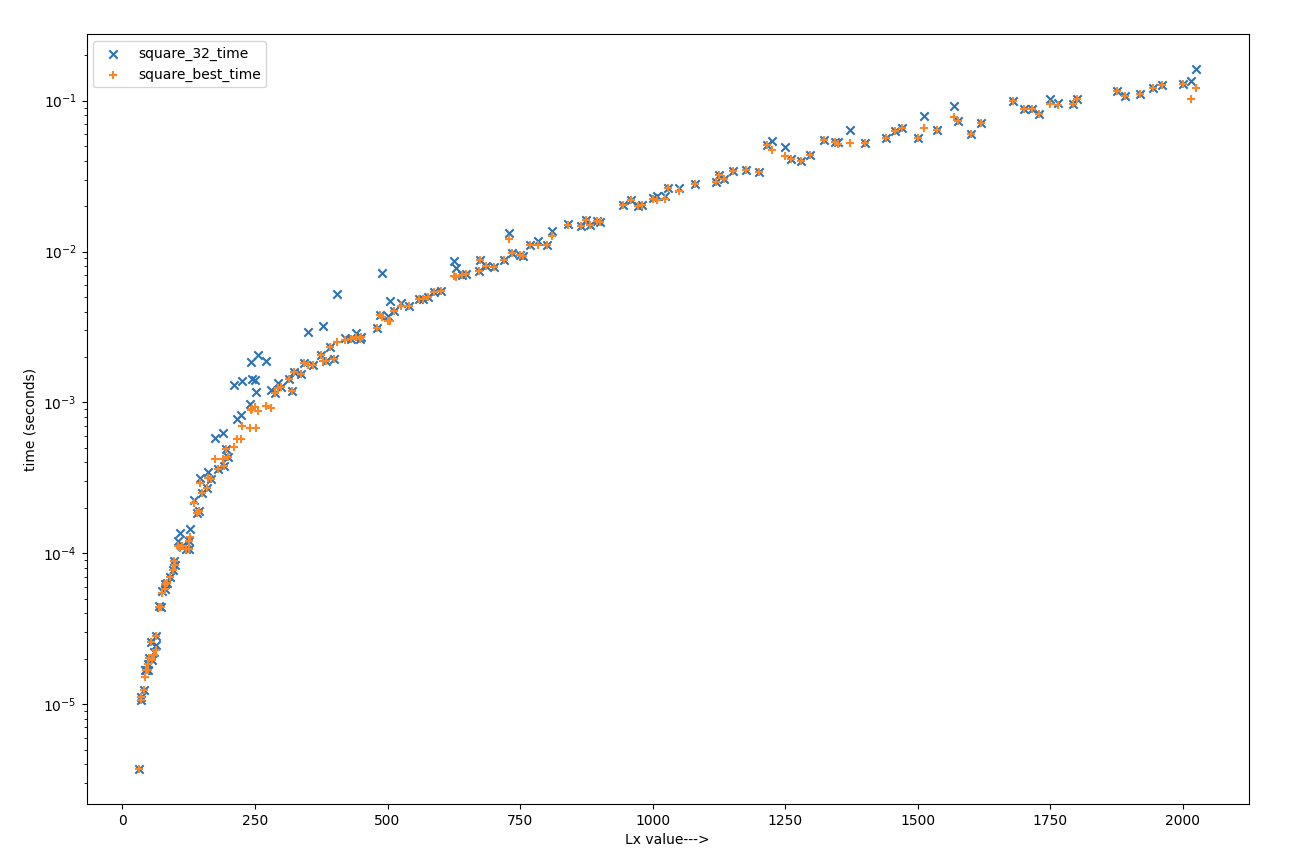} }}%
    \caption{Comparison of the skinny 32 rectangle vs the square parameters}%
    \label{fig:plotsfinal}%
\end{figure}

In addition to the observations made in the previous paragraph, it is worth noting that Figure \ref{fig:plotsfinal} shows that by using only the 32-rectangle parameters, we achieve timings that are extremely close to those of the square's best parameters. Moreover, this approach yields the expected $\mathcal{O}(L_x^2 \log L_x)$ time complexity for convolution. These findings lead us to propose that we can fix the height of our rectangles to 32 and use the best parameters to estimate extremely large squares while obtaining a satisfactory approximation. However, this is not the only height value we will explore. We plan to search for other height values, such as 4, 8, and 16, to determine if there are any other significant results.

\section{Future Work}
Although much of the foundational work has been completed over the past year, much work remains to be done. Some of them include 
\begin{itemize}
    \item Port to C++ and include in convolution libraries such as FFTW++; create wrappers for Python and Julia.
\item Create specialized algorithms for the hermitian and real case to save additional memory. This includes expanding the hybrid dealiasing codebase to support hermitian and real case convolutions.
\item Add code to generalize to $n$ dimensions.
\item In-depth study of generalizing the heuristics to problems apart from multiples of 2, 3, 5, 7.
\item Develop more efficient machine learning algorithms and study the classification algorithms in depth.
\end{itemize}

\section{Applications}

Our proposed research has many applications, especially since convolutions of different sizes are essential because they allow a model to detect patterns at different scales. For example, a small convolution filter may detect edges or minor details in an image. In contrast, a larger convolution filter may be able to detect overall shapes or structures. There are many potential applications of our hybrid dealiasing framework. One example is medical imaging, where accurate and efficient image reconstruction is essential for diagnosis and treatment planning. Our approach could improve the quality of medical images while reducing computation time and memory usage, leading to faster and more accurate diagnoses. Another potential application is in computer vision, where image processing is essential for tasks such as object recognition and tracking. Our framework could improve the accuracy of image processing algorithms while reducing computation time, leading to faster and more efficient computer vision systems. Lastly, we expect hybrid dealiasing to become the standard method for convolutions in the future. We hope to develop hybrid dealiased convolutions for computing sparse FFTs, fractional phase, and partial FFTs.

\section{Discussions}
This paper presented our work on expanding the Hybrid Dealiasing framework to address issues with uneven input or minimum padded sizes, reducing memory and computation time. We demonstrated the effectiveness of our approach in dealing with the unequal case in 1/2/3 dimensions and showed that our hybrid dealiasing solution is much quicker than typical explicit padding approaches. Furthermore, we extended our approach to multi-convolution and developed a conventional grid search hyperparameters search strategy that selects the optimal hyperparameters resulting in the shortest convolution time. We also developed efficient heuristics that allow us to estimate optimal parameters for larger convolution problems using only small squares/rectangles.

\justify
To improve the accuracy and efficiency of our approach, we developed a more optimized and efficient hyperparameter tuning algorithm called \textit{experience} that searched over a broader set of parameters to determine the fastest algorithm for a given convolution problem empirically. The results of our research show that our approach significantly reduces memory and computation time compared to existing state-of-the-art fast Fourier transform libraries such as FFTW++. In addition, our approach also offers a low-memory FFT algorithm that is faster.

\justify
In conclusion, our work presents a systematic framework that offers unique solutions to dealiasing issues. We have demonstrated the effectiveness of our approach in dealing with uneven input or minimum padded sizes, multi-convolution, and hyperparameter tuning. Our approach has the potential to significantly reduce memory and computation time for a wide range of convolution problems. It can be applied to various fields, including image processing, signal processing, and computational fluid dynamics. Future research directions include applying our approach to larger datasets and optimizing our heuristics for more efficient parameter estimation.

\section*{Acknowledgements}
I want to express my gratitude to Professor John Bowman and Noel Murasko for their guidance and supervision this semester. With their support, I was able to achieve what I did. I am also deeply thankful to Professor Nicolas Guay for his valuable feedback and for teaching such an excellent course. Finally, I am grateful to my peers and family for their unwavering support and encouragement throughout this challenging time.

\justify
Finally, I want to express my heartfelt thanks to everyone who has been a part of my life and helped me navigate difficult moments. I feel blessed to have the strength and wisdom to live meaningfully, and I am grateful to God for keeping me alive and guiding me along the way.

~\\

\pagebreak



\end{document}